\documentclass[final]{elsart}

\usepackage{amssymb,amsmath,setspace}
\usepackage[english,francais]{babel}

\let\-=\mathbf
\let\~=\mathrm
\let\<=\langle
\let\>=\rangle
\def\N{{\mathbb N}}

\newtheorem{theorem}{Theorem}[section]
\newtheorem{lemma}[theorem]{Lemma}
\newtheorem{proposition}[theorem]{Proposition}
\newtheorem{corollary}[theorem]{Corollary}

\renewcommand{\theequation}{\arabic{equation}}
\def\og{\leavevmode\raise.3ex\hbox{$\scriptscriptstyle\langle\!\langle$~}}
\def\fg{\leavevmode\raise.3ex\hbox{~$\!\scriptscriptstyle\,\rangle\!\rangle$}}

\journal{arxiv}
\begin{document}
\centerline{}
\begin{frontmatter}


\selectlanguage{english}
\title{A note on the Karhunen-Lo\`eve expansions for 
infinite-dimensional Bayesian inverse problems}


\selectlanguage{english}
\author[authorlabel1]{Jinglai Li}
\ead{jinglaili@sjtu.edu.cn}

\address[authorlabel1]{Institute of Natural Sciences, Department of Mathematics, and MOE Key Laboratory of Scientific and Engineering Computing, Shanghai Jiao Tong University,
Shanghai 200240, China}


\medskip
\begin{center}
\end{center}

\begin{abstract}
In this note, we consider the truncated Karhunen-Lo\`eve expansion for approximating solutions to infinite dimensional inverse problems. 
We show that, under certain conditions, the bound of the error between a solution and its finite-dimensional approximation can be estimated without the knowledge of the solution.





\end{abstract}
\end{frontmatter}


\selectlanguage{english}
\section{Introduction}
Nonparametric inverse problems have applications in many scientific or engineering problems, ranging from geophysical tomography~\cite{bui2013computational} 
to medical imaging~\cite{hanke2003recent}. 
In such problems the unknown that we want to determine is of infinite-dimension, for example, a function of space or time.

Identifying the unknown is usually cast as an optimization problem that needs to be solved numerically. 
Infinite-dimensional problems can not be solved directly with standard numerical techniques. A common practice is to first approximate the unknown with a finite-dimensional parameter, and then solve the resulting finite-dimensional problem numerically.
In particular, when the inverse problem is treated in a Bayesian framework, the Karhunen-Lo\`eve~(K-L) expansion~(\cite{papoulis2002probability}, Chapter 11) can be used to construct such a finite-dimensional approximation.
In the K-L method, the unknown is represented by a finite expansion of the eigenfunctions of the covariance operator of the prior measure. 

The K-L method has been long used to reduce the dimensionality in practical problems~\cite{mclaughlin1996reassessment,li2006efficient,marzouk2009dimensionality}; however, the use of it is never rigorously justified to the best of my knowledge.
To be specific,  it is unclear whether a fixed-dimensional representation can well approximate the solutions of the problem. 
In this note, we address the problem by proving that, if $u$ is a solution to the inverse problem defined as a minimizer to Eq~\eqref{e:minI}, 
 the error bound between $u$ and its finite K-L approximation can be estimated without the knowledge of $u$.

\section{Problem setup}
We consider the inverse problems in a Bayesian framework~(see \cite{stuart2010inverse} for a comprehensive overview of the Bayesian methods for infinite-dimensional inverse problems).
We assume the state space $X$ is a separable {Hilbert} space with inner product $\<\cdot,\cdot\>_X$.
 Our goal is to estimate $u\in X$ from some data $y$. 
The Bayes' formula in this setting should be interpreted as providing the Radon-Nikodym derivative
between the posterior measure $\mu$ and the prior measure $\mu_0$~\cite{Cotter2009,dashti_map_2013}:
\begin{equation}
\frac{d\mu}{d\mu_0}(u) = \exp(-\Phi(u)),
\end{equation}
where $\exp(-\Phi(u))$ is the likelihood function.
A typical example is to assume that the unknown $u$ is mapped to the data $y$ via a forward model  
$y = G(u)+\zeta$,
where $G:X\rightarrow R^d$ and $\zeta$ is a $d$-dimensional Gaussian noise with mean zero and covariance $C$. 
In this case $\Phi(u) = |C^{-\frac12}(Gu-y)|^2_2$. 

Next we assume a Gaussian prior is used. Namely we let $\mu_0$ be a zero-mean Gaussian measure defined on $X$ with covariance operator $Q$. 
Note that ${Q}$ is symmetric positive and of trace class.
$E=Q^{\frac12}(X)$ is a Hilbert space
with inner product 
\[\<\cdot,\cdot\>_E = \<{Q}^{-\frac12}\cdot,{Q}^{-\frac12}\cdot\>_X,\]
 which is known as the Cameron-Martin space associated with measure $\mu_0$.
Often we are only interested in a point estimate of $u$, rather than the posterior measure $\mu$ itself. 
To this end, as is shown in \cite{Cotter2009,dashti_map_2013}, the maximum a posterior (MAP) estimator of $u$ can be defined as the minimizers of the Onsager-Machlup functional over $E$:
\begin{equation}
\min_{u\in E}I(u) := \Phi(u) + \|u\|^2_E, \label{e:minI}
\end{equation}
where $\|u\|^2_E=\<u,u\>_E$.
Note that Eq.~\eqref{e:minI} can also be understood as a classic inverse problem where the cost function $\Phi(\cdot)$ is minimized 
with a Tikhonov regularization in the Hilbert space $E$~\cite{bissantz2004consistency}.

\section{Karhunen-Lo\`eve  representation}
Note that solving Eq.~\eqref{e:minI} directly involves inverting the operator $Q$, which can be rather challenging in practice. 
Alternatively, one can use substitution $u = Q^\frac12 x$ and rewrite Eq.~\eqref{e:minI} as
\begin{equation}
\min_{x\in X}J(x) := \Phi(Q^{\frac12}x) + \|x\|^2_X. \label{e:minJ}
\end{equation}
The following proposition states the equivalence of the two optimization problems.
\begin{proposition}
If $x$ minimizes $J(x)$ over $X$, $u=Q^\frac12x$ minimizes $I(u)$ over $E$,
and if $u$ minimizes $I(u)$ over $E$, $x = Q^{-\frac12}u \in X$ minimizes $J(x)$ over $X$.
\end{proposition}

\noindent\textit{Proof.}
We prove the proposition by contradiction. 
First it is easy to verify that, for any $x\in X$ and $u\in E$ satisfying $u=Q^{\frac12}x$, we have $I(u)=J(x)$.
Let $x$ be a minimizer $J(\cdot)$ over $X$, and assume $u=Q^\frac12x$ is not a minimizer of $I(\cdot)$ over $E$.
Namely, there exists an $u'\in E$ such that $I(u')<I(u)$. It follows directly that
$x'=Q^{-\frac12}u' \in X$
and  $J(x') < J(x)$, which contradicts that $x$ is a minimizer of $J$ over $X$.
Thus we have proved the first part of the proposition. The second part can be proved by following the same argument.\qed

Now we introduce the K-L expansion to reduce the dimensionality of Eq.~\eqref{e:minJ}.
We start with the following lemma~(\cite{da2006introduction}, Chapter~1):
\begin{lemma}
There exists a complete orthonormal \label{lm:eigens}
basis $\{e_k\}_{k\in\N}$ on $X$ and a sequence of non-negative numbers $\{\lambda_k\}_{k\in\N}$
such that ${Q} e_k = \lambda_k e_k$ and $\sum_{k=1}^\infty \lambda_k <\infty$, i.e., 
 $\{e_k\}_{k\in\N}$ and $\{\lambda_k\}_{k\in\N}$ being the eigenfunctions and eigenvalues of $Q$ respectively.
\end{lemma}

The basic idea of the K-L method is to solve the optimization problem in a finite-dimensional subspace of $X$: 
\begin{equation}
\min_{x\in X_n}J(x) := \Phi(Q^{\frac12}x) + \|x\|^2_X, \label{e:minJn}
\end{equation}
where $X_n$ be the space spanned by $\{e_k\}_{k=1}^n$ for a given $n\in\N$. 
In numerical implementation Eq.~\eqref{e:minJn} can be recast as
\begin{equation}
\min_{(\xi_1,\,...,\,\xi_n)\in R^n} \Phi(\sum_{k=1}^n \xi_k\sqrt{\lambda_k} e_k) + \sum_{k=1}^n \xi_k^2
\end{equation}
which is the usual K-L representation.
As is mentioned earlier, a critical question here is whether the finite subspace $X_n$ can provide good approximation to the solutions of Eq.~\eqref{e:minJ}.
Our main results regarding this problem are presented in the following theorem:
\begin{theorem} \label{th:representer}
Suppose $\Phi(u)$ is locally Lipschitz continuous, i.e., for every $r>0$, there exists a constant $L_r>0$ such
that for all $z_1,\,z_2\in X$ with $\|z_1\|_X,\,\|z_2\|_X<r$, we have
\[|\Phi(z_1)-\Phi(z_2)|< L_r \|z_1-z_2\|_X.\]
Let $\{e_k\}_{k\in\N}$ and $\{\lambda_k\}_{k\in\N}$ be the eigenfunctions and eigenvalues of $Q$ as defined in Lemma~\ref{lm:eigens}.
 There exists a constant $L>0$ such that, for any $x\in \arg\min_{x\in X} J(x)$, we have 
\[\|x-x_n\|_X<L\sqrt{\lambda^*_n}\] where 
$x_n = \sum_{k=1}^n \<x,e_k\>_X\, e_k$,  
and  $\lambda^*_n = \max_{k>n}\lambda_k$.

\end{theorem}

\noindent\textit{Proof.}
Let $x\in X$ be a minimizer of Eq.~\eqref{e:minJ}.
Since $\{e_k\}$ is a complete orthonormal basis for $X$, $x$ can be written as 
\[
x = \sum_{k=1}^\infty \xi_k\, e_k,  
\]
where $\xi_k= \<x,e_k\>_X$.
Let 
\[x_n = \sum_{k=1}^n \xi_k\, e_k. \]
As $x$ is a minimizer of $J(\cdot)$, take $r = \Phi(0)+1$ and so we have 
$J(x)< r$, which implies that $\|x_n\|_X\leq\|x\|_X<r$.
 $Q^\frac12$ is bounded, and so we have 
$\|Q^\frac12x_n\|_X,\,\|Q^\frac12x\|_X<\|Q^\frac12\| r$. 
Now recall that $\Phi(\cdot)$ is locally Lipschitz continuous, and so there exists a constant $L>0$ such
that 
\[
|\Phi(Q^\frac12x)-\Phi(Q^\frac12x_n)|< L \|Q^\frac12x-Q^\frac12 x_n\|_X.\]
Since $x$ minimizes $J(\cdot)$, we have $J(x)\leq J(x_n)$ which implies
\begin{align*}
\|x-x_n\|_X^2 &\leq |\Phi(Q^\frac12x)-\Phi(Q^\frac12x_n)|
< L \|Q^{\frac12}x-Q^{\frac12}x_n\|_X\\
&=L \<x-x_n,Q(x-x_n)\>^\frac12_X
=L\<\sum_{k=n+1}^\infty \xi_k e_k,\sum_{k=n+1}^\infty \xi_k \lambda_k e_k\>^{\frac12}_X
\\&= L \sqrt{\sum_{k=n+1}^\infty \lambda_k \xi_k^2}
\leq L\sqrt{\lambda^*_n}\|x-x_n\|_X.
\end{align*}
It then follows immediately that 
\[\|x-x_n\|_X\leq L\sqrt{\lambda^*_n}.\]\qed

Certainly we also want to know if the minimizer of the original problem~\eqref{e:minI} is well approximated by the K-L expansion.
To this end, we have the following corollary, which is a direct consequence of Theorem~\ref{th:representer}:
\begin{corollary}
Let $u= Q^\frac12 x$ and $u_n = Q^\frac12 x_n$, and we have $\|u-u_n\|_X < L \lambda^*_n$.
\end{corollary}

Another important question is that whether a solution to finite-dimensional problem~\eqref{e:minJn} is a good approximation to
that of the infinite-dimensional problem~\eqref{e:minJ}.
We have the following results regarding this issue:
\begin{corollary}
Let $x_n' \in \arg\min_{x\in X_n}J(x)$
and we have \[\min_{x\in X}J(x)\leq J(x_n')\leq \min_{x\in X}J(x)+L^2\lambda^*_n\]
\end{corollary}
The corollary follows directly from Theorem~\ref{th:representer} and so proof is omitted.

%
\section{Concluding remarks}
We theoretically study the truncated K-L expansions for approximating the solutions of infinite-dimensional 
Bayesian inverse problems. 
We show that the error between a solution to the inverse problem and its projection on the chosen finite-dimensional space 
is bounded by the eigenvalues of the covariance operator of the prior. 

\section*{Acknowledgment}
The work is supported by the NSFC under grant number 11301337.
\bibliographystyle{plain}
\bibliography{representer}

\begin{thebibliography}{10}

\bibitem{bissantz2004consistency}
N.~Bissantz, T.~Hohage, and A.~Munk.
\newblock Consistency and rates of convergence of nonlinear tikhonov
  regularization with random noise.
\newblock {\em Inverse Problems}, 20(6):1773, 2004.

\bibitem{bui2013computational}
T.~Bui-Thanh, O.~Ghattas, J.~Martin, and G.~Stadler.
\newblock A computational framework for infinite-dimensional bayesian inverse
  problems part i: The linearized case, with application to global seismic
  inversion.
\newblock {\em SIAM Journal on Scientific Computing}, 35(6):A2494--A2523, 2013.

\bibitem{Cotter2009}
S.~L. Cotter, M.~Dashti, J.~C. Robinson, and A.~M. Stuart.
\newblock {Bayesian inverse problems for functions and applications to fluid
  mechanics}.
\newblock {\em Inverse Problems}, 25(11):115008, November 2009.

\bibitem{da2006introduction}
G.~Da~Prato.
\newblock {\em An introduction to infinite-dimensional analysis}.
\newblock Springer, 2006.

\bibitem{dashti_map_2013}
M.~Dashti, K.~J.~H. Law, A.~M. Stuart, and J.~Voss.
\newblock {MAP} estimators and their consistency in {Bayesian} nonparametric
  inverse problems.
\newblock {\em Inverse Problems}, 29(9):095017, September 2013.

\bibitem{hanke2003recent}
M.~Hanke and M.~Br{\"u}hl.
\newblock Recent progress in electrical impedance tomography.
\newblock {\em Inverse Problems}, 19(6):S65, 2003.

\bibitem{li2006efficient}
W.~Li and O.~A Cirpka.
\newblock Efficient geostatistical inverse methods for structured and
  unstructured grids.
\newblock {\em Water resources research}, 42(6), 2006.

\bibitem{marzouk2009dimensionality}
Y.~M Marzouk and H.~N Najm.
\newblock Dimensionality reduction and polynomial chaos acceleration of
  bayesian inference in inverse problems.
\newblock {\em Journal of Computational Physics}, 228(6):1862--1902, 2009.

\bibitem{mclaughlin1996reassessment}
D.~McLaughlin and L.~R. Townley.
\newblock A reassessment of the groundwater inverse problem.
\newblock {\em Water Resources Research}, 32(5):1131--1161, 1996.

\bibitem{papoulis2002probability}
A.~Papoulis and S.~U. Pillai.
\newblock {\em Probability, random variables, and stochastic processes}.
\newblock Tata McGraw-Hill Education, 2002.

\bibitem{stuart2010inverse}
A.~M. Stuart.
\newblock Inverse problems: a {Bayesian} perspective.
\newblock {\em Acta Numerica}, 19:451--559, 2010.

\end{thebibliography}
\end{document}